\documentclass[11pt]{article}

\usepackage{amssymb}
\usepackage{amsthm}
\usepackage{amsmath}
\usepackage{amsopn}
\usepackage{amstext}
\usepackage{epsfig}

\DeclareMathOperator{\interior}{Int}

\DeclareMathOperator{\dist}{dist}

\def\beq{\begin{equation}}
\def\eeq{\end{equation}}

\def\Om{\ensuremath{\Omega} }
\def\Dom{\ensuremath{\partial\Omega} }

\def\om{\omega}

\def\D{\ensuremath{\mathbb D} }
\def\R{\ensuremath{\mathbb R} }
\def\C{\ensuremath{\mathbb C} }

\def\N{\ensuremath{\mathbb N} }
\def\P{\ensuremath{\mathbb P} }
\def\Z{\ensuremath{\mathbb Z} }

\def\P{\ensuremath{\mathbb P}}

\theoremstyle{plain}
\newtheorem{thm}{Theorem}[section]
\newtheorem{lem}[thm]{Lemma}

\newtheorem{cor}[thm]{Corollary}
\newtheorem{prop}[thm]{Proposition}

\theoremstyle{definition}
\newtheorem{defn}{Definition}[section]

\theoremstyle{remark}

\numberwithin{equation}{section}

\newcommand{\lemref}[1]{Lemma~\ref{#1}}

\newcommand{\ccl}[1]{{\mathbf {#1}}} 

\newcommand{\cl}{\operatorname{cl}}

\newcommand{\length}{\operatorname{length}}

\renewcommand{\Re}{\operatorname{Re}}
\renewcommand{\Im}{\operatorname{Im}}

\newcommand{\PP}{{\mathbf P}}
\newcommand{\EE}{{\mathbf E}}
\newcommand{\CC}{{\mathbb C}}

\begin{document}

\title{On the computational complexity of the Riemann mapping}         

\author{I.~Binder\thanks{partially supported by  NSERC Discovery grant}\\Dept. of Mathematics\\ University of Toronto  \and M.~Braverman\thanks{partially supported by NSERC 
Post-graduate Scholarship}\\Dept. of Computer Science\\ University of Toronto \and M.~Yampolsky
\thanks{partially supported by  NSERC Discovery grant}\\Dept. of Mathematics\\ University of Toronto}
\date{\today}          

\maketitle

\begin{abstract}
In this paper we consider the computational complexity 
of uniformizing a domain with a given computable boundary. 
We give nontrivial upper and lower bounds in two settings:
when the approximation of the boundary is given either as a list of pixels, or 
 by a Turing Machine.
\end{abstract}


\section{Introduction}

\subsection{Foreword.} Computational conformal mapping 
is prominently featured in problems of applied analysis and
mathematical physics, as well as in engineering disciplines,
such as image processing. In this paper we address the theoretical
foundations of numerically approximating the conformal mapping
between two planar domains. We obtain a lower bound on the
computational complexity of an algorithm solving this problem,
and show that this bound is almost sharp. To achieve the
latter, we present a very space-efficient probabilistic algorithm for
constructing such a mapping.

\subsection*{Acknowledgments} The second author would like 
to thank Ker-I Ko for bringing the problem to 
his attention during CCA'04. 

The authors are grateful to Stephen Cook for his numerous 
helpful suggestions during the preparation of the paper. 

\subsection{Background in computational complexity theory}

We present here some basic definitions and results from the 
computational complexity theory. A more comprehensive
discussion can be found in \cite{Sip,Papa}.

The primary goal of the computational complexity theory
is to classify different computational problems into 
{\em complexity classes} according to their computational 
hardness. The basic abstract object here is a {\em Turing 
Machine} which for most purposes can be thought of as a program in any 
programming language.

The complexity class $\ccl{P}$ includes problems that are
computable in time polynomial in the length of the input. 
Those are thought of as the ``relatively easy" problems. 
Examples of problems in $\ccl{P}$ include arithmetic operations, 
finding a shortest path in a graph and primality testing. 
``Difficult" problems, such as factoring integers or computing the optimal strategy for 
playing ``Go" on an $n\times n$ board, are generally
thought not to be in $\ccl{P}$. Whether a problem is in $\ccl{P}$ 
or not is usually a good criterion in assessing its
true hardness. By an analogy with $\ccl{P}$ one 
can define the class $\ccl{EXP}$ of problems solvable in time 
$2^{n^c}$ for some $c$ on input length $n$. Factoring integers is
in $\ccl{EXP}$ via an obvious exhaustive search. Playing ``Go" optimally is  also in $\ccl{EXP}$, since we can easily 
enumerate all possible games and compute an optimal path in time $2^{O(n^2)}$.
Using a diagonalization argument, it is not hard to see that $\ccl{P} \subsetneq  \ccl{EXP}$ (see e.g. \cite{Sip}).

The complexity class $\ccl{NP}$ contains problems that 
are easy to verify, but may be hard to guess. More precisely, 
a predicate $Q(x)$ is in $\ccl{NP}$, if there is a poly-time 
computable predicate $R(x,y)$, where $y$ has length
polynomial in the length of $x$, such that $Q(x) = \exists y~R(x,y)$.
By an exhaustive search for $y$ one sees that $\ccl{NP} \subseteq \ccl{EXP}$.
There is a subclass of $\ccl{NP}$ called the $\ccl{NP}$-complete
problems, or $\ccl{NPC}$. Problems in $\ccl{NPC}$ have the property 
of being the ``hardest" in $\ccl{NP}$: if one could solve {\em any} 
problem in $\ccl{NPC}$ in polynomial time, then one could solve {\em all}
$\ccl{NP}$ problems in polynomial time. 

One of the most famous   $\ccl{NP}$-complete problems is the 
{\em satisfiability} problem $\text{\it SAT}$. The problem is the 
following: given a propositional formula $\phi(y)$ does 
it have a truth assignment $y_T$ such that $\phi(y_T)=1$. 
An example of a problem in $\ccl{NP}$ that is thought to be 
hard but not $\ccl{NP}$-complete is the following. Given a 
pair of numbers $m<n$, determine whether $n$ has a divisor 
between $2$ and $m$. This problem can be used to factor integers. 
It is in $\ccl{NP}$ since it can be formulated as $\exists k~(1<k<m)\wedge
k \mid n$. It is one of the Clay \$1,000,000 questions of whether 
$\ccl{P}=\ccl{NP}$. 

A bigger class of problems is the class $\ccl{\#P}$. It is the class
of problems which are equivalent to counting the number of 
satisfying assignments for a given propositional formula -- 
this natural complete problem for this class is denoted by $\#\text{\it SAT}$.
Obviously $\ccl{NP} \subseteq \ccl{\#P}$, since to solve $\text{\it SAT}$ we only need
to know whether the number of its satisfying assignments is 
bigger than $0$ or not, which is easier than actually determining 
this number.

The next class of problems is the class $\ccl{PSPACE}$ -- the class 
of problems solvable in {\em space} polynomial in the input size. 
It is easy to see that all the classes mentioned above are in
$\ccl{PSPACE}$. On the other hand, $\ccl{PSPACE} \subseteq \ccl{EXP}$, since a machine 
with $p(n)$ memory bits can have at most $2^{p(n)}$ different configurations, 
and can run for at most $2^{p(n)}$ steps without getting into an infinite 
loop.

The class of problems solvable in {\em logarithmic} space is a 
class of problems that are solvable in space $O(\log n)$ for 
input size $n$. Here the input and the output are read-only and 
write-only respectively. This class is denoted by $\ccl{L}$. By the same 
reasoning as $\ccl{PSPACE} \subseteq \ccl{EXP}$, we have $\ccl{L}\subseteq \ccl{P}$. A {\em randomized} 
version of $\ccl{L}$ are the problems that can be solved correctly 
with error probability $<1/n$ in space $O(\log n)$ and time 
$\operatorname{poly}(n)$. This class is called $\ccl{BPL}$.   
 
Overall, we have the following chain of inclusions:
$$
\ccl{L} \subseteq \ccl{BPL} \subseteq \ccl{P} \subseteq  \ccl{NP} \subseteq \ccl{\#P} \subseteq \ccl{PSPACE} \subseteq \ccl{EXP}.
$$
By diagonalization, $\ccl{BPL} \neq \ccl{PSPACE}$, and $\ccl{P} \neq \ccl{EXP}$. No other separations 
are known. 

In recent years, some progress has been made in {\em derandomizing} 
$\ccl{BPL}$, that is in showing that there is a deterministic 
algorithm that does not require a lot of additional computational resources. 
We will need the following recent result on
 derandomization:

\begin{thm} \cite{nisan}
\label{derandomization}
 There exists a deterministic
algorithm for the following problem:

\medskip
\noindent
{\bf Input:} An $n \times n$ transition probability matrix $M$, 
an integer $t$, and a rational $\epsilon$. 

\medskip
\noindent
{\bf Output:} A matrix $A$ such that $|| A - M^t || \leq \epsilon$.

The algorithm runs in time $\text{poly}(N)$ and space $O(\log^2 N)$, where 
$N = n + t + \epsilon^{-1}$.
\end{thm}

We will also use a {\em circuit complexity} class. A circuit 
consist of inputs, logical gates, and an output. The gates are usually 
NOT (one input, one output), AND, and OR. The latter gates can either 
have two, or unboundedly many inputs. In the discussion below, any number 
of inputs is allowed. The {\em size} of a circuit is the number of gates used.
The {\em depth} of a circuit is the number of gates on the longest path 
from an input to the output. It is known that any boolean function $f: \{0,1\}^n
\rightarrow \{0,1\}$ can be computed by a circuit of size $O(2^n/n)$.
It is not hard to see that functions in $\ccl{P}$ are computable by polynomial 
size circuits. The class $\ccl{AC^0}$ is the class of functions that 
are computable by a family of circuits (one for each input size), 
that have   constant depth and polynomial size. This is one of the very few
complexity classes for which non-diagonalization lower bounds exist. 
In particular it has been shown that computing 
the parity of the number of $1$'s in a string cannot be done in $\ccl{AC^0}$ (see, for example, \cite{FSS}). 
A more general problem that cannot be done in $\ccl{AC^0}$ is the {\em majority problem} $\text{MAJ}_n$:
given a string $x\in \{0,1\}^n$, $\text{MAJ}_n(x)$ is $1$ if and only if the majority 
of the entries in $x$ are $1$.

\subsection{Computational complexity of sets}

We review the definition and the basic properties 
of computable sets. We refer the reader to \cite{BratWei,Wei,RetWeiSTOC,BravFOCS} for
a more comprehensive exposition. 

Intuitively, we say the time complexity of 
a set $S$ is $t(n)$ if it takes time $t(n)$ to decide whether to draw a pixel of size $2^{-n}$ 
in the picture of $S$.  Mathematically, the definition is as follows:

\begin{defn}
\label{def1}
A set $T$ is said to be a $2^{-n}$-picture of a bounded set $S$ if:
\begin{itemize}
\item[(i)] $S \subset T$, and 
\item[(ii)] $T \subset B(S,2^{-n}) =
\{ x \in \R^2~:~|x-s|<2^{-n}~\mbox{for some~} s\in S\}$.
\end{itemize}
\end{defn}

\noindent
Definition \ref{def1} means that $T$ is a $2^{-n}$-approximation of $S$ with respect to 
the {\em Hausdorff metric}, given by
$$
d_H (S,T) := \inf \{ r : S \subset B(T,r)~~\mbox{and}~~
T \subset B(S,r) \}.
$$

Suppose we are trying to generate a  picture of a set $S$ using 
a union of round pixels of radius $2^{-n}$ with centers at all 
the points of the form $\left( \frac{i}{2^{n}}, 
\frac{j}{2^{n}}\right)$, with $i$ and $j$ integers. In order to 
draw the picture, we have to decide for each pair $(i,j)$ whether to 
draw the pixel centered at $\left( \frac{i}{2^{n}},
\frac{j}{2^{n}}\right)$ or not. We want to draw the pixel if 
it intersects $S$ and to omit it if some neighborhood of the pixel 
does not intersect $S$. Formally, we want to compute a function 
\beq
\label{star1}
f_S ( n, i/2^{n}, j/2^{n}) = 
\left\{ 
\begin{array}{ll}
1, & B((i/2^{n}, j/2^{n}), 2^{-n})\cap S \neq \emptyset \\
0, & B((i/2^{n}, j/2^{n}), 2 \cdot 2^{-n})\cap S = \emptyset  \\
0~~\mbox{or}~~1,~~& \mbox{in all other cases}
\end{array}
\right. 
\eeq


The time complexity of  $S$ is defined as follows.

\begin{defn}
\label{defcomp}
A bounded set $S$ is said to be computable in time $t(n)$ if there 
is 
a function $f(n,\bullet)$ satisfying \eqref{star1} which runs in time 
$t(n)$. We say that $S$ is poly-time computable if there is a polynomial 
$p$, such that $S$ is computable in time $p(n)$. 
\end{defn} 

Computability of sets in bounded space is defined in a 
similar manner. There, the amount of {\em memory} the 
machine is allowed to use is restricted. 

To see why this is the ``right" definition, suppose we are trying to 
draw a set $S$ on a computer screen which has a $1000\times 1000$ pixel
resolution. A $2^{-n}$-zoomed in picture of $S$ has $O(2^{2 n})$ pixels
of size $2^{-n}$, and thus would take time $O(t(n) \cdot 2^{2 n})$ to 
compute. This quantity is exponential in $n$, even if $t(n)$ is bounded 
by a polynomial. But we are drawing $S$ on a finite-resolution display, 
and we will only need to draw $1000 \cdot 1000 = 10^6$ pixels. Hence 
the running time would be $O(10^6 \cdot t(n)) = O(t(n))$. This running 
time is polynomial in $n$ if and only if $t(n)$ is polynomial. 
Hence $t(n)$ reflects the `true' cost of zooming in.

\subsection{Background in complex analysis}
\label{subsec:icomp}

To make the paper self-contained, we list here a few results from complex analysis that will be used later. We refer to \cite{Ahlfors,Duren,pombook} for a more comprehensive discussion.
 
Let $\Om\subsetneq\CC$ be a simply-connected planar domain with $w\in\Om$. {\it The Riemann Uniformization Theorem} states that there is unique conformal map $\psi$ of $\Om$  onto the unit disk $\D$ with $\psi(w)=0$, $\psi'(w)>0$. The number $r(\Om,w)=1/\psi'(w)$ is called the {\it conformal radius} of $\Om$. Roughly speaking, $r(\Om,w)$ measures the size of $\Om$ as viewed from $w$:

\begin{prop}[Koebe's Theorem] 
In this notation we have
$$\dist(w,\partial\Om)\geq\frac{r(\Om,w)}{4}.$$
\end{prop}

We note the following basic monotonicity 
property of the conformal radius: 
\begin{lem}
\label{radius-monotone}
If $\Om_1\subset\Om_2$, $w\in\Om_1$, then $r(\Om_1,w)\leq r(\Om_2,w)$.
\end{lem}

By a theorem of Carath{\'e}odory (see e.g. \cite{pombook}),
if the boundary $\partial \Omega$ is a Jordan curve, then the map $\psi$ can be extended to a homeomorphism between the
 closure of $\Om$ and the closed unit disk $\cl(\D)$.

Let $z^{*}=1/\overline{z}$ be the inversion of $z$ with respect to the unit circle $\{|z|=1\}$. 
We will make use of the following particular case of the {\it Reflection Principle}: 

\begin{lem}
\label{reflection}
If $J\subset\{|z|=1\}$ is an open arc, and $\phi$ is a continuous map on  $\D\cup J$ which is analytic on $\D$, and $\phi(J)\subset\{|z|=1\}$, then the map $\Phi$ defined by 
\beq
\Phi(z)=
\begin{cases}
\phi(z),& |z|<1\text{ and }z\in J\\
\phi^{*}(z^{*}),&|z|>1
\end{cases} 
\eeq
is analytic on the domain $\D\cup\{|z|>1\}\cup J$.
\end{lem}

\noindent
In particular, if $\phi$ is a conformal map of \D onto a domain $\Om\subset\D$ with Jordan boundary, and $K$ is an open arc, $K\subset\partial\Om\cap\{|z|=1\}$, then $\Phi$ is a conformal map of $\D\cup\{|z|>1\}\cup J$, where $J=\phi^{-1}(K)$ ($\phi$ is extendable to $\cl(\D)$ by Carath{\'e}odory theorem).

Let now $\Omega$ be a domain with the boundary $\Dom$ consisting of finitely many Jordan curves. Let $f$ be a continuous function on $\Dom$. 
A function $u:\cl(\Omega)\to\CC$ is a {\it solution for the Dirichlet problem with the boundary data $f$},
if
\begin{itemize}
\item $u$ is continuous in $\cl(\Om)$, 
\item $u$ is harmonic on $\Om$ ($\Delta u=\partial_{xx}u+\partial_{yy}u=0$), and 
\item $u(z)=f(z)$ for $z\in\Dom$. 
\end{itemize}

\noindent
For any $f$ such a solution exists and is unique.
Moreover, there exists a unique family of  measures $\om_{w,\Om}$ on $\Dom$ such that for any $f\in C(\Dom)$, 
$$u(w)=\int_{\Dom}f(z)\,d\om_{w,\Om}(z).$$
The measure $\om_{w,\Om}$ is called a {\it harmonic measure}.
If one fixes $K\subset\Dom$, the function $w\mapsto\om_{w,\Om}(K)$ is harmonic in $\Om$.

\noindent
If $\Om$ is simply-connected, then for a set $K\subset\Dom$, we have 
$$\om_{w,\Om}(K)=\frac{1}{2\pi}\length(\psi(K)),$$ where $\psi$ is the Riemann map of $\Om$ onto $\D$ with $\psi(w)=0$.

The Dirichlet problem can be solved probabilistically. Namely, for $w\in\Om$, let $B_w(t)$ be the two-dimensional Brownian motion started at $w$, and the exit time be defined by $$T=\inf\{t\ :\ B_w(t)\not\in\Om\}.$$ 
Then the solution of the Dirichlet problem is given by the following formula of Kakutani (see e.g. \cite{GM}): 
$$u(w)=\EE[f(B_w(T))].$$
Note that the harmonic measure for a set $K\subset\Dom$ is now given by
$$\om_{w,\Om}(K)=\PP[B_w(T)\in K].$$

We will make use of the {\it Maximum Principle for harmonic functions} (see \cite{Ahlfors}): 
\begin{lem}
\label{maximum principle}
If $u_1(z)$ and $u_2(z)$ are two functions which are harmonic in $\Om$, continuous on the whole $\cl(\Om)$, and  $u_1(z)\geq u_2(z)$ for $z\in\Dom$, then $u_1(z)\geq u_2(z)$ for all $z\in\Om$.
\end{lem}

\noindent
An easy consequence is the {\it Monotonicity Property of the harmonic measure} (see \cite{pombook}): 
\begin{cor}
\label{monotonicity}
If $w\in\Om_1\subset\Om_2$, $K\subset\Dom_1\cap\Dom_2$, then $$\om_{w,\Om_1}(K)\leq\om_{w,\Om_2}(K).$$
\end{cor}

\noindent
We will also make use of a Distortion Theorem for conformal maps(see \cite{Duren}):

\begin{thm}
\label{distortion}
If $\phi$ is conformal in the disk $\{z\ :\ |z-w|<r\}$, then
\beq\label{eq:distfunc}
|\phi'(w)|\frac{r^2|z-w|}{(r+|z-w|)^2}\leq |\phi(z)-\phi(w)|\leq |\phi'(w)|\frac{r^2|z-w|}{(r-|z-w|)^2}
\eeq
and
\beq\label{eq:distder}
r^2|\phi'(w)|\frac{r-|z-w|}{(r+|z-w|)^3}\leq |\phi'(z)|\leq r^2|\phi'(w)|\frac{r+|z-w|}{(r-|z-w|)^3}\eeq
\end{thm}

\subsection{Results}

In Section \ref{sec:up} we propose a new algorithm for computing 
the Riemann map. We use the random walks solution to the general Dirichlet
problem to produce a solution to the uniformization problem. This 
gives an extremely space-efficient algorithm. 

The formulation of the theorem will depend on how the boundary of the uniformized domain 
$\Omega$ is specified for our algorithm. Since the domain $\Omega$ we consider is computable,
there exists a Turing machine $M(n)$ which for a given $n$ computes a function \eqref{star1}.
Our algorithm may then query $M(n)$ for different values of $(i/2^n,j/2^n)$ to ascertain
whether this particular dyadic rational point lies within one-pixel distance from $\partial \Omega$.
A formal way of saying this is that our algorithm will have an access to an {\it oracle} for a function
given by \eqref{star1}.

\begin{thm}
\label{thm:ub}
There is an algorithm $A$ that computes the uniformizing map in 
the following sense. 

Let  $\Om$ be a bounded simply-connected domain, and $w_0 \in \Om$. 
$\Dom$ is provided to $A$ by an oracle representing it in 
the sense of equation \eqref{star1}. Then $A$ computes the absolute values of the uniformizing map 
$\phi: (\Om, w_0) \rightarrow (\D, 0)$
 with precision $2^{-n}$ in space bounded by $C \cdot n^{2}$, and time
 $2^{O(n)}$, where $C$ depends only on the diameter of
$\Om$ and $d(w_0, \partial \Om)$. Furthermore, the algorithm computes the value of $\phi(w)$ with precision $2^{-n}$ as long as $|\phi(w)|<1-2^{-n}$. Moreover, $A$ queries $\Dom$
with precision of at most $2^{-O(n)}$.

In particular, if $\partial\Om$ is
 polynomial space computable in space $n^a$ for some constant 
$a\geq1$ and time $T(n)<2^{O(n^a)}$, then $A$ can be used to compute
the uniformizing map in space $C \cdot n^{\max(a,2)}$ and 
time $2^{O(n^a)}$. 
\end{thm}

In the scale where the entire boundary is given to us explicitly, and 
not by an oracle for it, we have the following.
 
\begin{thm}
\label{thm:ubsm}
There is an algorithm $A'$ that computes the uniformizing map in 
the following sense. 

Let  $\Om$ be a bounded simply-connected domain, and $w_0 \in \Om$. 
Suppose that for some $n=2^k$, $\partial\Om$ is given to $A'$ 
with precision $\frac{1}{n}$ by $O(n^2)$ pixels.
 Then $A'$ computes the absolute values of the uniformizing map $\phi: (\Om, w_0) \rightarrow (\D, 0)$
 within an error of $O(1/n)$ in randomized 
space  bounded by $O(k)$ and time polynomial in $n=2^k$ (that is, by a $\ccl{BPL}(n)$-machine). Furthermore, the algorithm computes the value of $\phi(w)$ with precision $1/n$ as long as $|\phi(w)|<1-1/n$.
\end{thm}

In Section \ref{sec:low}, we show that even if the domain we are uniformizing 
is very simple computationally, the complexity of the uniformization
can be quite high. Moreover, it might be difficult to compute the 
conformal radius of the domain.

More specifically, the following theorems are established in the Section \ref{sec:low}.
\begin{thm}
\label{thm:lb}
Suppose there is an algorithm $A$ that given a simply-connected domain
$\Om$ with a  linear-time computable boundary and an inner radius $>\frac{1}{2}$ and a number $n$ computes
the first $20n$ digits of the
conformal radius $r(\Om,0)$, then we can use one call to $A$ to solve any
instance of a $\#\text{\it SAT}(n)$ with a linear time overhead. 

In other words, $\ccl{\#P}$ is poly-time reducible to computing the conformal
radius of a set. 
\end{thm}

\begin{thm}
\label{thm:lbsm}
Consider the problem of computing the conformal radius
 of a simply-connected domain $\Om$, where the 
boundary of $\Om$
 is given with precision $\frac{1}{n}$
by an explicit collection of $O(n^2)$ pixels.

Denote the problem of computing the conformal radius with precision 
$\frac{1}{n^{c}}$ by $\text{CONF}(n, n^c)$.
 Then $\text{MAJ}_n$ is $\ccl{AC^0}$ 
reducible to $\text{CONF}(n, n^c)$ for any $0<c<1/2$.
\end{thm}

\subsection{Comparison with known results.}

\medskip
The first constructive proof of the Riemann Uniformization Theorem is due to Koebe \cite{Koebe},
and dates to the early 1900's. Formal proofs of the constructive nature of the Theorem
which follow Koebe's argument 
under various computability conditions on the boundary of the domain are numerous in the
literature (see e.g. \cite{Cheng,BishBr,Zhou,Hertling}).
In particular, Zhou \cite{Zhou} and Hertling \cite{Hertling} give constructive proofs
under computability conditions on the boundary similar to those used by us.
The question of complexity bounds on the construction was raised, in particular, in most
of the works quoted above. However, the only result known to us was announced by Chou in 
\cite{Chou}. He states that in the case when the boundary is poly$(n)$ computable,
the problem of computation of the mapping is in $\ccl{EXPSPACE}(n)$.

\medskip

From the practical (that is, applied) point of view, the most computationally
efficient algorithm used 
nowadays to calculate the conformal map is the ``Zipper'', invented by
Marshall (see \cite{Mar}). 
The effectiveness of this algorithm was recently studied by Marshall and Rohde in \cite{MR}. 
The ``Zipper'', however, falls beyond the theoretical upper bound on the complexity
of this problem, which we establish in  Section \ref{sec:up}: 
in the settings of the 
Theorem \ref{thm:ub}, it computes the  uniformizing map in space 
$2^{O(n^a)}$ and time $2^{O(n^a)}$, and thus belongs to the complexity class $\ccl{EXP}$.
It is reasonable to expect then, that an algorithm can be found in the class $\ccl{PSPACE}$
which is more practically efficient than ``Zipper''.


\section{Computing the uniformization in polynomial space}
\label{sec:up}
Let $\Om$ be a bounded simply-connected planar domain, let $K\subset\Omega$ be a fixed compact set with smooth boundary with $\dist(K, \partial\Omega)>10 \cdot 2^{-n}$. First we discuss a probabilistic algorithm for solving the Dirichlet problem in the domain $\Omega\setminus K$ with precision $2^{-n}$.

\subsection{General Dirichlet problem}

The discrete analogue of the Dirichlet problem can be defined as follows. 
For $H\subset h\Z^2$ ($h>0$), the interior of $H$ is defined by $\interior(H)=\{a\in H\ :\ a\pm h, a\pm ih\in H\}$. The boundary of $H$ is defined by $\partial H=h\Z^2\setminus \left(\interior(H)\cup \interior(h\Z^2\setminus H)\right)$.
We say that a function $u$ defined on  $H\subset h\Z^2$ is {\it discrete harmonic} if for any $a\in \interior(H)$ we have 
$$u(a)=1/4(u(a+h)+u(a-h)+u(a+ih)+u(a-ih)).$$

\noindent
Let $B^w_n$ be the standard Random Walk (cf \cite{Spitzer})
on $h\Z^2$ started at $w\in H$, where $H$ is closed ($\partial H\subset H$). Let the exit time $N$ be defined as $N=\min\{n\ :\ B^w_n\not\in H \}-1$.
Let $f$ be a function on $\partial H$.
It is almost obvious that the function
$$u(w)=\EE[f(B^w_N)]$$
is discrete harmonic on $H$. This function is called the {\it solution of the Dirichlet problem with boundary data $f$} (cf. the continuous case discussion in subsection \ref{subsec:icomp}).

Let $\Om$ be a domain with boundary $\Dom$ and $f\in C(\Dom)$. For $h>0$ define $H_h=\Om\cap h\Z^2$. For $w\in\partial H_h$, let $f_h(w)=f(z)$, where $z$ is one of the points on $\Dom$ closest to $w$. The solution $u_h$ of the corresponding discrete Dirichlet problem, is called the $h$-discrete solution of the original continuous Dirichlet problem.

We need the following easy case of the approximating property of the $h$-discrete solutions (see \cite{Spitzer}, \cite{Laa}).

\begin{lem}
\label{lem:discdir}
Let $\Om$ be a domain. Let $f$ be a continuous and locally constant function on $\Dom$ taking only $0$ and $1$ values, and $u$ is the solution of the corresponding Dirichlet problem. Let $h<10^{-3}$ be such that for any $z_1, z_2\in\Dom$ with $|z_1-z_2|<\sqrt{h}$ we have $f(z_1)=f(z_2)$. Let $u_h$ be the $h$-discrete solution. Then if $\dist(w, \Dom)>\sqrt{h}$, then $|u(w)-u_h(w)|\leq 2\sqrt{h}$.
\end{lem}

Since the exit probabilities of a random walk can be computed by a $\ccl{BPL}(h^{-1})$ 
machine; 
if the values of $f$ and the boundary
 $\partial H_h$ are given by an oracle, then $u$ can be computed in the randomized space 
$O(-\log h)$ and time $O(h^{-2})$.
Thus Lemma \ref{lem:discdir} immediately implies the following statement about the solution 
of the general Dirichlet problem:

\begin{lem}
\label{GenDir}
There is a randomized algorithm $D$ that computes a solution of the Dirichlet
problem in the following sense. 

Let $\Om$ be a bounded planar domain and $K\subset\Omega$ 
be a fixed compact set with smooth boundary and $\dist(K, \Dom)>10 \cdot 2^{-n}$.
Suppose that $f$ is the function which is equal to $0$ on $\Dom$ and 
$1$ on $K$. Then $D$ computes the solution of the corresponding Dirichlet
problem with precision $2^{-n}$, $2^{-n}$-away
from $\partial \Omega \cup K$ in space $O(n)$, and time $2^{O(n)}$. 
The computation is done probabilistically, 
and outputs the correct value within an error of $2^{-n}$ with probability $>\frac{1}{2}$. 

In particular,
if both $K$ and $\Dom$ are computable in space $n^a$ for some 
constant $a\geq1$ and time $T(n)<2^{O(n^a)}$.
Then we can compute the solution of the Dirichlet problem for any point,
 which is at least $2^{-n}$ away from $\partial \Om$ and $K$
in space $O(n^a)$, and time $2^{O(n)} T(n)$. 
\end{lem}

\subsection{The conformal radius}
Let  $w_0\in\Omega$, and let $\psi$ be the conformal mapping of $\Om$ onto the unit disk $\D$ with $\psi(w_0)=0$ and $\psi'(w_0)>0$.
Assume that $\partial\Om$ is given to us within an error of $2^{-n}$ in Hausdorff metric and that $d(w_0, \partial \Om)\ge 1$.
As a first application of Lemma \ref{GenDir} let us give an algorithm for calculating $|\psi'(w_0)|$ with precision $2^{-n}$ in space $O(n^a)$, and time $2^{O(n)} T(n)$. 
Denote $$w_1 = w_0 + e^{-n}\text{ and }K_1 = B(w_0, e^{-2n})$$ 

\begin{lem}\label{ConfRad}
Let $h_1$ be the solution of the following Dirichlet problem: 
$$
\begin{cases}
h_1(w)=1,&|w-w_0|=e^{-2n}\\
h_1(w)=0,&w\in\Dom\\
\Delta h_1(w)=0,&w\in\Om\setminus K_1
\end{cases}
$$
Then 
$$\left|\log|\psi'(w_0)|-n\left(\frac{1-2h_1(w_1)}{1-h_1(w_1)}\right)\right|\leq 5n\ e^{-n}.$$
\end{lem}

\begin{proof}

By the first statement of Theorem \ref{distortion},
\beq\label{eq:include}
B(0, e^{-2n}/(1+e^{-2n})^2\psi'(w_0))\subset\psi(K_1)\subset B(0, e^{-2n}/(1+e^{-2n})^2\psi'(w_0))
\eeq 
Let $B_1=\psi^{-1}\left(B(0, e^{-2n}(1-3e^{-2n})\psi'(w_0))\right)$ and  $B_2=\psi^{-1}\left(B(0, e^{-2n}(1-3e^{-2n})\psi'(w_0))\right)$.
Since $1/(1+e^{-2n})^2<(1-3e^{-2n})$ and $(1+3e^{-2n})<1/(1-e^{-2n})^2$, \eqref{eq:include} implies
\beq\label{eq:dist}
B_1\subset K_1\subset B_2
\eeq  

The functions $$H_1(w)=\frac{\log|\psi(w)|}{-2n+\log(1-3e^{-2n})+\log\psi'(w_0)}$$ and $$H_2(w)=\frac{\log|\psi(w)|}{-2n+\log(1+3e^{-2n})+\log\psi'(w_0)}$$ are harmonic in $\Om\setminus B_1$ and $\Om\setminus B_2$ respectively, equal to $0$ on $\Dom$, and equal to $1$ on the boundaries of $B_1$ and $B_2$ respectively.
By the Maximum Principle, $H_1\leq h_1 \leq H_2$, or, more explicitly,
\beq\label{eq:max}
\frac{\log|\psi(w)|}{-2n+\log(1-3e^{-2n})+\log\psi'(w_0)}\leq h_1(w)\leq \frac{\log|\psi(w)|}{-2n+\log(1+3e^{-2n})+\log\psi'(w_0)}
\eeq

Another application of the same distortion theorem yields
\beq\label{eq:distest}
e^{-n}(1-3e^{-n})\psi'(w_0)\leq|\psi(w_1)|\leq e^{-n}(1+3e^{-n})\psi'(w_0).
\eeq
Evaluating both sides of the inequality \eqref{eq:max} at the point $w_1$ using \eqref{eq:distest} completes the proof of the lemma.
\end{proof}

It now follows from Lemma \ref{GenDir} that we can compute $\psi'(w_0)$ with the same complexity constraints as in Lemma \ref{GenDir}. 

\subsection{The Riemann map}

Let $h_1$, $K_1$ be as in the previous section.
\begin{lem}\label{Dist}\label{AbsValue}
Let $|w-w_0|>e^{-n}$.
Then
$$\left|\log|\psi(w)|-h_1(w)(\log\psi'(w_0)-2n)\right|\leq 3\cdot e^{-2n}.$$
\end{lem}
\begin{proof}
By the equation \eqref{eq:max}, 
$$h_1(w)\log(1-3e^{-2n})\leq\log|\psi(w)|-h_1(w)(\log\psi'(w_0)-2n)\leq h_1(w)\log(1+3e^{-2n})$$
To prove the lemma it suffices to notice that $h_1(w)\leq 1$ and $|\log(1+x)|\leq |x|$.
\end{proof}

\noindent
Using Lemmas \ref{ConfRad} and \ref{GenDir}, we see that $|\psi(w)|$ is computable with the same restrictions as in Lemma \ref{GenDir}, provided that $\dist(w,\Dom)>e^{-n}$ and $|w-w_0|>e^{-n}$.

Now we have to compute $\arg\bigl(\psi(w)\bigr)$. 
To achieve this, we introduce another Dirichlet problem. Let $K_2=B(w_0+e^{-2n}, e^{-4n})$, and 
let $h_2$ be the solution of the following Dirichlet problem: 
$$
\begin{cases}
h_2(w)=1,&|w-w_0-e^{-2n}|=e^{-4n}\\
h_2(w)=0,&w\in\Dom\\
\Delta h_2(w)=0,&w\in\Om\setminus K_2
\end{cases}
$$

\begin{figure}
\centerline{\includegraphics[angle=0,scale=0.6]{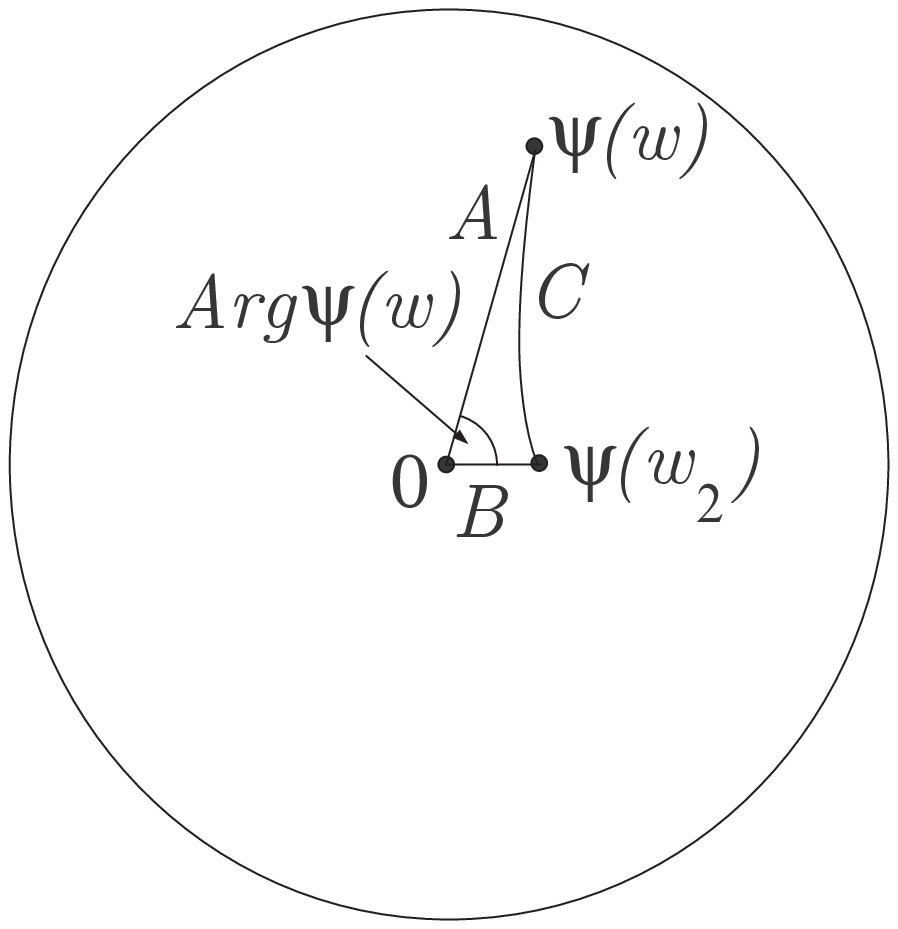}}
\caption{Computing $\arg\bigl(\psi(w)\bigr)$}
\label{Arg}\end{figure}

\noindent
Define $$\tilde\psi(w)=\frac{\psi(w)-e^{-2n}\psi'(w_0)}{1-\psi(w)e^{-2n}\psi'(w_0)}.$$ $\tilde\psi$ is also a Riemann map from $\Om$ onto $\D$. Let $w_2:=\tilde\psi^{-1}(0)=\psi^{-1}(e^{-2n}\psi'(w_0))$. By Distortion Theorem \ref{distortion}, $$|w_2-w_0-e^{-2n}|\leq 2\cdot e^{-4n}.$$ 
As in Lemma \ref{Dist},
$$\left|\log|\tilde\psi(w)|-h_2(w)(\log\tilde\psi'(w_2)-4n)\right|\leq 3\cdot e^{-4n}.$$
Another application of the Distortion Theorem \ref{distortion} yields $|\log|\tilde\psi'(w_2)|-\log|\tilde\psi'(w_0+e^{-2n})||<2e^{-4n}$. So, finally,
$$\left|\log|\tilde\psi(w)|-h_2(w)(\log|\tilde\psi'(w_0+e^{-2n})|-4n)\right|\leq 5\cdot e^{-4n}.$$
\noindent
Now we use a standard formula from hyperbolic trigonometry (see \cite{Thurston})
$$
\cos\arg\psi(w)=\frac{\cosh C -\cosh A\cosh B}{\sinh A\sinh B},
$$
where $$A=\log\frac{1+|\psi(w)|}{1-|\psi(w)|},\ B=\log\frac{1+e^{-2n}\psi'(0)}{1-e^{-2n}\psi'(0)}\text{, and }C=\log\frac{1+|\tilde\psi(w)|}{1-|\tilde\psi(w)|}.$$
See Figure \ref{Arg}.

Note that $\sinh B\sim e^{-2n}$, $\sinh A>e^{-n}$ when $|w|>e^{-n}$, $\cosh B-1\sim e^{-2n}$. Using the error estimate in the Lemma \ref{Dist}, we obtain that the formula allows us to compute $\cos\arg\psi(w)$ up to $e^{-n}$, provided that $|\psi(w)|<1-e^{-n}$.

Using the same argument for computing $\cos\arg(\psi(w)/i)$, we can completely determine the value of $\arg\phi(w)$.

Now we can give an algorithm which satisfies the conditions of Theorems \ref{thm:ub} and \ref{thm:ubsm}. 
\begin{proof}[Proof of Theorems \ref{thm:ub} and \ref{thm:ubsm}]
We can create a $\text{poly}(n)\times \text{poly}(n)$ matrix $M$ representing the transition
probabilities between the poly$(n)$ possible states of the random walk.
Simulating the random walk for $t=\text{poly}(n)$ steps amounts to approximating 
$M^t$. The required precision is also inverse polynomial in $n$. 
By Theorem \ref{derandomization}, this can be done in time polynomial 
in $n$, and space $O(\log^2 n)$, which imply Theorem \ref{thm:ubsm}. By changing the scale, and replacing 
$n$ with $2^n$, we obtain Theorem \ref{thm:ub}.
\end{proof}


\section{Lower bounds on the complexity of uniformization}\label{sec:low}

In this section we establish Theorems \ref{thm:lb} and \ref{thm:lbsm}.

Let us first remark that by Distortion Theorem \ref{distortion}, any algorithm computing values of the uniformization map will also compute the conformal radius with the same precision.
 
\begin{figure}
\begin{center}
\includegraphics[angle=0, scale=0.35]{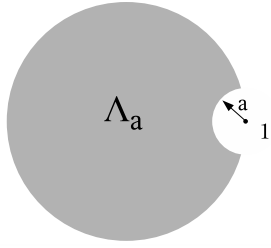}
\end{center}
\caption{$\Lambda_a$}
\label{LA}\end{figure}

Let $\Lambda_a$ be the domain $\D\setminus \{|z-1|\leq a\}$ -- the unit disk with a small bump of radius $a$ removed (see Figure \ref{LA}).

Fix a large $n\in\N$. Let now for $0\leq l<2^n$, and let $\Om_l=e^{2\pi i l/2^n}\Lambda_{2^{-10n}}$ be the rotated domain $\Lambda_{2^{-10n}}$. For a set $L=\left\{l_1, l_2, \dots, l_k\right\}$ with all $0\leq l_1<l_2\dots l_k<2^n$, let $\Om_L=\Om_{l_1}\cap\Om_{l_2}\cap\dots\Om_{l_k}$. Thus $\Om_L$ is the unit disk with $k$ relatively ``spread out'' bumps removed. 

\begin{thm}\label{thm:main}
For large enough $n$, 
$$\left|r(\Om_L, 0)-1+k2^{-20n-1}\right|<\frac1{10}2^{-20n}.$$
\end{thm}

To prove Theorem \ref{thm:main},
 we estimate the conformal radius of $\Lambda_a$ for an arbitrary $a$. 

\begin{lem}\label{lem:bump}
The conformal radius of $\Lambda_a$ is equal to $\displaystyle\frac{2-2a}{2-2a+a^2}.$
\end{lem}

\noindent
As a consequence we get that for large $n$,
\begin{equation}\label{eq:bumpn}
|r(\Lambda_{2^{-10n}})-1+2^{-20n-1}|<2^{-30n+2}
\end{equation}
\begin{proof}[Proof of Lemma \ref{lem:bump}]
Let $\P=\C\setminus\left\{\Im z=0,\ \Re z\leq0\right\}$ be the complex plane with the negative real axis removed.

The function $$\chi(z)=\left(\frac{1-z}{1+z}\right)^2$$ maps $\D$ conformally onto $\P$ , $\chi(0)=1$. It also maps $\Lambda_a$ onto $\Lambda'_b=\P\setminus\{|z|\leq b\}$, where $b=\displaystyle\left(\frac{a}{2-a}\right)^2$.

Observe also that  
$$h(z)=\frac{z+b^2/z-2b}{(1-b)^2}$$ 
maps $\Lambda'_b$ conformally onto $\P$, with $h(1)=1$, $\displaystyle h'(1)=\frac{1+b}{1-b}=\frac{2-2a+a^2}{2-2a}$.

Thus the map $\phi_0(z)=\chi^{-1}\circ h^{-1}\circ \chi(z)$ maps $\D$ conformally onto $\Lambda_a$, and the conformal radius of $\Lambda_a$ is equal to
$$r(\Lambda_a)=1/h'(1)=\frac{2-2a}{2-2a+a^2}.$$
\end{proof}

For a set $L$, let $\phi_L$ be the conformal map of $\D$ onto $\Om_L$ with $\phi_L(0)=0$, $\phi'_L(0)>0$ ($\phi_L$ is the inverse of the uniformization map).
Let $L'=(l_2, l_3,\dots, l_k)$ be the set $L$ with the first element removed.  Let $g(z)=\phi_{L'}^{-1}\circ\phi_L(z)$ be the conformal map of $\D$ onto 
$$\Gamma=\D\setminus\phi_{L'}^{-1}\left(\{|z|<1,\ |z-e^{2\pi i l_1/2^n}|\leq2^{-10n}\}\right).$$
Let us also introduce two domains $$\Gamma_{+}=\D\setminus\{z\ :\ |w-z|\leq2^{-10n}(1-2^{-2n})\}$$ and $$\Gamma_{-}=\D\setminus\{z\ :\ |w-z|\leq2^{-10n}(1+2^{-2n})\},$$ where $w=\phi_{L'}^{-1}(1)$.

We will use the following property of $\Gamma$
\begin{lem}\label{lem:distest}
$$\Gamma_{-}\subset\Gamma\subset\Gamma_{+}.$$
\end{lem}

\noindent
Let us first show how to derive Theorem \ref{thm:main} from Lemma \ref{lem:distest}.

\noindent
By Lemma \ref{lem:distest} and \lemref{radius-monotone} (monotonicity of conformal radius), 
$$r(\Gamma_{-})\leq g'(0)=r(\Gamma)\leq r(\Gamma_{+})$$

\noindent
Now Lemma \ref{lem:bump} implies that for large $n$
$$|r(\Gamma_{-})-1+2^{-20n-1}|<2^{-22n+2},\quad |r(\Gamma_{+})-1+2^{-20n-1}|<2^{-22n+2},$$
and thus 
$$|g'(0)-1+2^{-20n-1}|<2^{-22n+2}.$$

\noindent
Note now that $\phi_L(z)=\phi_{L'}\circ g(z)$, so $r(\Om_{L})=g'(0)r(\Om_{L'})$.
The Theorem easily follows from this relation by induction on the size of $L$.

\noindent
So to establish   Theorem \ref{thm:main}, it is enough to prove Lemma \ref{lem:distest}.

\begin{figure}[t]
\begin{center}
\includegraphics[angle=0, scale=0.65]{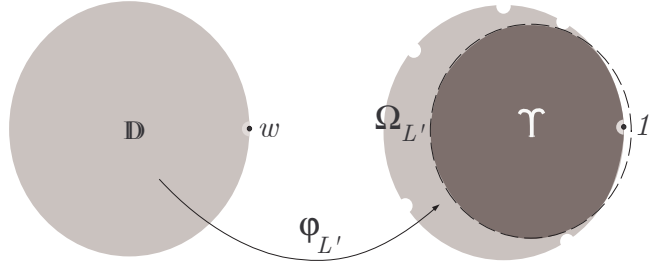}
\end{center}
\caption{The map $\phi_{L'}$ and domain $\Upsilon$}.
\label{Phi}\end{figure}

\begin{proof}[Proof of Lemma \ref{lem:distest}]

Without loss of generality we can assume that $l_1=0$.

Let $\Upsilon=\D\cap\{|z-2^{-7n+1}|<1-2^{-7n}\}$. Note that $\Upsilon\subset\Om_{L'}$, since $0\not\in L'$.
Let $\psi$ be the conformal map of $\D$ onto $\Upsilon$ with $\psi(0)=0$, $\psi'(0)>0$.

Let $K$ be the arc $[1,e^{\pi i 2^{-7n}}]$. Note that $K\subset\partial \Upsilon\cap\partial \Om_{L'}$. Let $K'=\phi_{L'}^{-1}(K)$. The normalized length of $K'$, $|K'|$, is the harmonic measure of $K$ in $\Om_{L'}$ evaluated at zero. By monotonicity of harmonic measure (\lemref{monotonicity}) it is bounded above by $2^{-7n-1}=|K|$ and below by the harmonic measure of $K$ in $\Upsilon$ evaluated at zero.

So
\beq\label{eq:hmeas}
2^{-7n-1}\geq|K'|\geq|\psi^{-1}(K)|\geq 2^{-7n-1}(1-2^{-5n})
\eeq
The same estimate applies to $K''=\phi_{L'}^{-1}([e^{-\pi i2^{-7n}},1])$.

Using these estimates we see that the arc $$ [e^{-\pi i2^{-7n}(1-2^{-5n})},e^{\pi i2^{-7n}(1-2^{-5n})}])$$ is mapped by $\phi_{L'}$ inside $K'\cup K''\subset\partial\D$. Thus we can use the Reflection Principle \ref{reflection} to extend  the map $\phi_{L'}$ to a map $G$ of the whole disk  $\{|z-w|<2^{-7n-1}(1-2^{-5n})\}$ with $G(w)=1$.

Let $J_{\epsilon}$ be the arc $[e^{-i\epsilon},e^{i\epsilon}]$. 
Using the fact that $\Upsilon\subset\Om_{L'}$ and the monotonicity of harmonic measure, we obtain
\beq
|G^{-1}(J_{\epsilon})|=|\phi_{L'}^{-1}(J_{\epsilon})|=\om_{0,\Om_{L'}}(J_{\epsilon})> \om_{0,\Upsilon}(J_{\epsilon})=|\psi^{-1}(J_{\epsilon})|\geq |J_{\epsilon}|(1-2^{-5n})
\eeq
Letting $\epsilon\to0$, we obtain
\beq\label{eq:bounddist}
|G'(w)|=\lim_{\epsilon\to0}\frac{|J_{\epsilon}|}{|G^{-1}(J_{\epsilon})|}>1-2^{-5n}
\eeq

Now we can use the Distortion Theorem \ref{distortion} applied to the disk $\{|z-w|<2^{-7n}(1-2^{-5n})\}$ to see that
$$G(\{z\ :\ |w-z|\leq2^{-10n}(1-2^{-2n})\})\subset\{|z-1|<2^{-10n}\}$$
and
$$\{|z-1|<2^{-10n}\}\subset G(\{z\ :\ |w-z|\leq2^{-10n}(1+2^{-2n})\})$$

\begin{figure}
\begin{center}
\includegraphics[angle=0, scale=0.4]{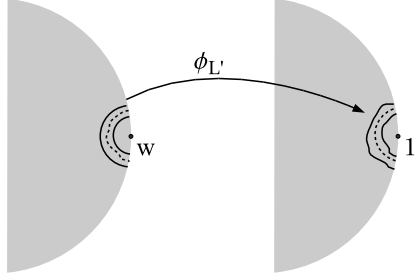}
\end{center}
\caption{$\phi_{L'}$ in the neighborhood of $w$}
\label{magn}\end{figure}
But this is precisely the statement of the lemma.
\end{proof} 
Now we are in the position to prove Theorems \ref{thm:lb} and \ref{thm:lbsm}.
\begin{proof}[Proof of Theorem \ref{thm:lb}]
For a propositional formula $\Phi$ with $n$ variables, let $L\subset \{0, 1, \ldots, 2^n-1\}$ be the set of 
numbers  corresponding to
its satisfying instances. Then the boundary of $\Om_L$ is computable in 
linear time, given the access to $\Phi$.  Theorem \ref{thm:main} now  implies that using 
$r(\Om_L, 0)$ we can evaluate $|L|=k$, and solve the $\#SAT$ problem on $\Phi$,
which is exactly Theorem \ref{thm:lb}. 
\end{proof}

\begin{proof}[Proof of Theorem \ref{thm:lbsm}]
Suppose that we are given a string $s$ of $n=2^k$ zeros and ones. We can view 
it as a set $L \subset \{0, 1, \ldots, 2^k-1\}$. $\Om_L$ can be obtained from 
$L$ by a trivial one-layered circuit with just NOT gates. Theorem \ref{thm:main}
implies that using $r(\Om_L, 0)$ with $2^{-O(k)}$ precision, we can evaluate
$|L|$ and solve the MAJ$_n$ problem on $s$,
which is exactly Theorem \ref{thm:lbsm}. 
\end{proof}

\bibliographystyle{amsalpha}

\end{document}